\documentclass[12pt,oneside, italian]{article}

\usepackage[english]{babel}

\usepackage{latexsym}

\usepackage{amssymb,amsthm,amsmath}

\usepackage{graphicx}

\title{\bf The Spectrum of Quantum Caps in $PG(4,4)$}

\author{Daniele Bartoli, Stefano Marcugini, Fernanda Pambianco}

\newtheorem{theorem}{Theorem}[section]

{\theoremstyle{definition}

\newtheorem*{definition*}{Definition}

\newtheorem{definition}[theorem]{Definition}}

\newtheorem*{proposition*}{Proposition}

\newtheorem*{corollary*}{Corollary}

\newtheorem*{lemma*}{Lemma}

\date{}

\begin{document}

\maketitle \vspace*{-15mm} \noindent

\ \\

\ \\

\begin{abstract}

\noindent We prove the non existence of quantum caps of sizes $37$ and $39$. This completes the spectrum of quantum caps in $PG(4,4)$. This also implies the non existence of linear $[[37,27,4]]$ and $[[39,29,4]]$-codes. The problem of the existence of non linear quantum codes with such parameters remains still open.

\end{abstract}

\section{Introduction}

In the projective space $PG(r,q)$ over the Galois Field $GF(q)$, an $n$-cap is a set of $n$ points no $3$
of which are collinear. An $n$-cap is called \emph{complete} if it is not contained in an $(n+1)$-cap. \\
We call an $n$-cap an $n$-\emph{quantum} cap if the code generated by its matrix is a
\emph{quantum stabilizer code}, that is:
\begin{definition}\label{def:quantumcodes}
A \emph{quaternary quantum  stabilizer code} is an additive quaternary code $\mathcal{C}$
contained in its dual $\mathcal{C}^{\bot}$, where the duality is with respect to the symplectic form.
\end{definition}
\noindent In particular:
\begin{definition}\label{parametri}
A quantum code  $\mathcal{C}$ with parameters $n,k,d$ ($[[n,k,d]]$-code), where $k>0$,
is a quaternary quantum stabilizer code of binary dimension $n-k$ satisfying the following:
any codeword of $\mathcal{C}^{\bot}$ having weight at most $d-1$ is in $\mathcal{C}$.\\
The code is \emph{pure} if $\mathcal{C}^{\bot}$ does not contain codewords of weight $< d$,
equivalently if $\mathcal{C}^{\bot}$ has strength $t \geq d-1$.\\
An $[[n,0,d]]$-code $\mathcal{C}$ is a self-dual quaternary quantum stabilizer code of strength
$t=d-1$.
\end{definition}
\noindent For a more detailed introduction of quantum codes see in particular \cite{tesi}, \cite{AMS}, \cite{NewQCaps}.\\
In $1999$ Bierbrauer and Edel (see \cite{41Cap}) have proven that the maximum size of complete caps in $PG(4,4)$ is $41$ and only two non equivalent $41$-caps exist. One of them results to be quantic.\\
In $2003$ (see \cite{40Cap}) they also have proven that in $AG(4,4)$ there exists a unique $40$-cap. It results to be quantic in $PG(4,4)$.\\
In $2008$ Tonchev has constructed quantum caps of sizes $10,12,14-27,29,31,33,35$ (see \cite{Tonchev08}), starting from the complete $41$-quantum cap in $PG(4,4)$.\\
In $2009$ we have found examples of quantum caps of sizes $13,28,30,32,34,36,38$, see \cite{AMS} and \cite{NewQCaps}. Moreover we have proven that there are only two examples of non equivalent quantum caps of size $10$ and five of size $12$. In the same article we have proven by exhaustive search that no $11$-quantum cap exists.\\
In this article we show that no quantum caps of sizes $37$ and $39$ exist. Therefore the following theorem holds:
\begin{theorem}
If $\mathcal{K} \subset PG(4,4)$ is a quantum cap, then $10\leq |\mathcal{K}| \leq 41$, with $|\mathcal{K}| \neq 11,37,39$.
\end{theorem}

\section{The searching algorithm }\label{alg}
We performed an exhaustive search for quantum caps of sizes $37$ and $39$. To do this, programs in C/C++ have been utilized.\\
We start from caps, complete and incomplete, in $PG(3,4)$ where the classification is known (see \cite{BMP-PG(3-4)-rep} and \cite{FMMP-ArsCombin}) and we try to extend every starting cap joining new points in $PG(4,4)$, to obtain complete or incomplete caps of sizes $37$ and $39$. We utilize the following geometric characterization to reduce the number of cases to examine (\cite{NewQCaps}, Theorem $3.4$). 
\begin{theorem}\label{equivalenza}  The following are equivalent:
\begin{enumerate}
\item A pure quantum $[[n,k,d]]$-code which is linear over $\mathbb{F}_{4}$.
\item A set of $n$ points in $PG(\frac{n-k}{2}-1,4)$ of strength $t=d-1$, such that the intersection
size with any hyperplane has the same parity as $n$.
\item An $[n,k]_{4}$ linear code of strength $t=d-1$, all of whose weights are even.
\item An $[n,k]_{4}$ linear code of strength $t=d-1$ which is self-orthogonal with respect to the
Hermitian form.
\end{enumerate}
\end{theorem}
According to the previous theorem we can consider starting caps in $PG(3,4)$ of odd size only.\\
In particular we consider in our search only caps of sizes $13$, $15$ and $17$ in $PG(3,4)$, since the following theorem and the non existence of particular linear codes.
\begin{theorem}
The following are equivalent:
\begin{enumerate}
\item An $[n,k,d']_{q}$-code with $d' \geq d$.
\item A multiset $\mathcal{M}$ of $n$ points of the projective space $PG(k-1,q)$, satisfying the following: for every hyperplane $H \subset PG(k-1,q)$
there are at least $d$ points of $\mathcal{M}$ outside $H$ (in the multiset sense).
\end{enumerate}
\end{theorem}
More precisely, we know that linear codes with $n =37,39$ $k=5$ and $d > n-12$ do not exist
(see \cite{codetable}) and so there exists an hyperplane which contains at least 12 points of the caps.\\
Then we consider only the examples of non equivalent caps in $PG(3,4)$ contained in the following table:\\
\begin{table}[h]
\caption{Number and type of non equivalent caps $\mathcal{K} \subset PG(3,4)$, with $|\mathcal{K}|=13,15,17$}
\begin{center}
\begin{tabular}{|c|c|c|}
\hline
$|\mathcal{K}|$&\# COMPLETE&\# INCOMPLETE\\
&CAPS&CAPS\\
\hline
13&1&3\\
15&0&1\\
17&1&0\\
\hline
\end{tabular}
\end{center}
\end{table}

\noindent We proceeded in this way:
\begin{enumerate}
\item we start from all non equivalent caps, complete and incomplete, in an hyperplane of $PG(4,4)$ and we extend them by the addition of new points of
$PG(4,4)\setminus PG(3,4)$;
\item we obtain the caps of sizes $37$ and $39$ in $PG(4, 4)$ by an exhaustive search;
\item we control if an obtained cap is a quantum cap, according to Theorem \ref{equivalenza}. In particular we check if all the weights of the linear code generated by the cap are even.
\end{enumerate}

\section{Results}\label{risultati}
We finish our search, finding no examples of quantum caps in $PG(4,4)$ of sizes $37$ and $39$. According \cite{AMS}, \cite{NewQCaps}, \cite{41Cap}, \cite{40Cap} and \cite{Tonchev08} we have proven the following:
\begin{theorem}
If $\mathcal{K} \subset PG(4,4)$ is a quantum cap, then $10\leq |\mathcal{K}| \leq 41$, with $|\mathcal{K}| \neq 11,37,39$.
\end{theorem}

\end{document}